\documentclass[preprint]{elsarticle}
\usepackage[latin1]{inputenc}
\usepackage[swedish,english]{babel}
\usepackage{amsmath,amsthm}
\usepackage{amssymb}
\usepackage{verbatim}
\usepackage{graphicx}
\usepackage[pdftitle={High order 2D panels},
  pdfauthor={Goude and Engblom},
  pdffitwindow=true,
  breaklinks=true,
  colorlinks=true,
  urlcolor=blue,
  linkcolor=red,
  citecolor=red,
  anchorcolor=red]{hyperref}

\numberwithin{equation}{section}
\numberwithin{table}{section}
\numberwithin{figure}{section}

\theoremstyle{plain}

\theoremstyle{definition}

\theoremstyle{remark}

\begin{document}

\begin{frontmatter}
\title{A general high order two-dimensional panel method}

\author[doe]{Anders Goude}
\ead{anders.goude@angstrom.uu.se}

\author[dosc]{Stefan Engblom\corref{cor1}}
\ead{stefane@it.uu.se}

\address[doe]{Division of Electricity, Department of Engineering
    Sciences, Uppsala University, SE-751~21 Uppsala, Sweden.}

\address[dosc]{Division of Scientific Computing, Department of
    Information Technology, Uppsala University, SE-751~05 Uppsala,
    Sweden. }

\cortext[cor1]{Corresponding author, telephone +46-18-471 27 54, fax +46-18-51 19 25.}
\date{\today}


\selectlanguage{english}

\begin{abstract}

  We develop an efficient and high order panel method with
  applications in airfoil design. Through the use of analytic work
  and careful considerations near singularities our approach is
  quadrature-free. The resulting method is examined with respect to
  accuracy and efficiency and we discuss the different trade-offs in
  approximation order and computational complexity. A reference
  implementation within a package for a two-dimensional fast multipole
  method is distributed freely.

\medskip


\medskip




\end{abstract}
\end{frontmatter}

\section{Introduction}
\label{sec:Introduction}

Most methods to numerically solve partial differential equations
(PDEs) fall into one of two categories. The first is volume
discretization methods, including, for example, finite element and
finite volume methods. Here the resulting set of equations is large
but sparse since the discretization nodes are connected only
locally. If the PDE has a known fundamental solution, the full
solution may instead be obtained in a discretization procedure
involving only the boundary. This approach is commonly known as
boundary element methods (BEMs) and generally involves fewer unknowns
which, however, are connected globally.

In the current work, focus will be on fluid mechanical applications
where the Laplace equation is used to calculate potential flow
solutions. For many aerodynamics simulations, the target is a small
object in the form of, e.g., an airfoil in a large domain which makes
the BEM particularly attractive \cite{Drela89, Bal1998343,
  maskew1982prediction}. For time-dependent calculations, the method
is commonly combined with the release and subsequent advection of
\emph{vortices}, effectively discretization points which approximates
the flow. This is the method known as the \emph{vortex method}
\cite{Cottet08}. Here the flow velocity $\vec{V}$ is calculated as a
combination of a potential flow and of vortex contributions,
\begin{equation}
  \label{eq:VMVelocity}
  \vec{V} = \nabla\phi+\vec{V}_{\omega},
\end{equation}
where $\phi$ is the solution to a Laplace equation and where
$\vec{V}_{\omega}$ is the contribution from the vortices in the flow
\cite{maskew1982prediction, Eldredge2007626, Zanon13}. Note that the
contribution to the flow velocity from each vortex has to be
calculated at each vortex position and at every time step. This is an
\emph{$N$-body problem} for which, as discussed below, fast algorithms
should be employed.

One common method for solving the potential flow problem is to
discretize the boundaries using \emph{panels}. For fluid mechanical
applications, these panels are constructed such that the flow satisfy
the no-penetration boundary condition,
\begin{equation}
  \label{eq:NoPenetration}
  \vec{V}\cdot\hat{n} = 0.
\end{equation}
The standard approach in 2D is to discretize the boundary using linear
panels consisting of straight line segments, and to use a panel
strength which is either constant or linearly varying along this
line. These panels can be constructed from source and vortex sheets,
but other possibilities also exist \cite{drela2014flight}.

An obvious issue with panels that have a linear shape is that there
will be sharp corners at the transitions between panels such that the
solution to Laplace's equation approaches infinity at the corner. For a pure
potential flow solution, this is often not a problem, since the no
penetration boundary condition is only satisfied at the centers of the
panels. For vortex methods, however, this yields large numerical
errors when evaluating the velocity in the vicinity of such transition
points. A remedy is to use panels with higher order shapes to make
the boundary smooth. For three-dimensional implementations,
such panels have been discussed in \cite{willis2006quadratic,
  Gao20081271, moore2013progress}. Here, numerical integration is
required to solve the flow contribution from the panel. This can be
time consuming, especially for velocity evaluations close to the
singularity of the fundamental solution, and it is therefore desirable
to avoid numerical integration whenever possible. Indeed, for
two-dimensional calculations, Ramachandran and co-authors have derived
a solution involving panels with cubic shape and a linear distribution
of the panel strength \cite{Ramachandran03}.

To improve the computational speed, vortex methods commonly rely on
the fast multipole method (FMM) \cite{Greengard87}, and the same
method can also be used to accelerate the solution of the dense BEM
matrix with the influence coefficients of the panels
\cite{Liu09,Ramachandran03}.

In the present paper we extend the work in \cite{Ramachandran03} and
design a general framework for two-dimensional panels with high order
in both shape and strength and which does not require numerical
integration. We develop the necessary analytic relations in
\S\ref{sec:theory}, where we also discuss practical implementation
issues allowing the method to be evaluated via the FMM. The
performance of our framework is evaluated in \S\ref{sec:evaluation}
and conclusions are summarized in \S\ref{sec:conclusions}.

Our method has been implemented within the 2D FMM-software described
in \cite{Engblom11, Goude13a}, and is distributed as open source. See
\S\ref{subsec:reproducibility} for details.


\section{Theory and implementation}
\label{sec:theory}

We construct the panels from point sources/vortices in
\S\ref{subsec:background}. The procedure for calculating the
contribution from a panel is developed in \S\ref{subsec:velocity},
including both near and far-field evaluations. How the boundary
conditions can be solved is described in \S\ref{subsec:bc}, which
includes how to integrate the contribution of a source point/panel
over the panel surface. Finally, corrections to ensure continuity of
the source strength are discussed in \S\ref{subsec:smooth}

\subsection{Representation of panels}
\label{subsec:background}

For a two dimensional flow using the complex number representation, the velocity $V$ at position $z$ from a
source/vortex at position $z_{v}$ is given by
\begin{equation}
  \label{eq:vortexvelocity}
  \overline{V\left(z\right)}=\frac{1}{2\pi}\frac{Q+i\Gamma}{z-z_{v}},
\end{equation}
where $\overline{V\left(z\right)}$ denotes the complex conjugate of
the velocity, $Q$ is the source strength, and $\Gamma$ is the vortex
strength.

\begin{figure}
  \includegraphics[width=0.8\columnwidth]{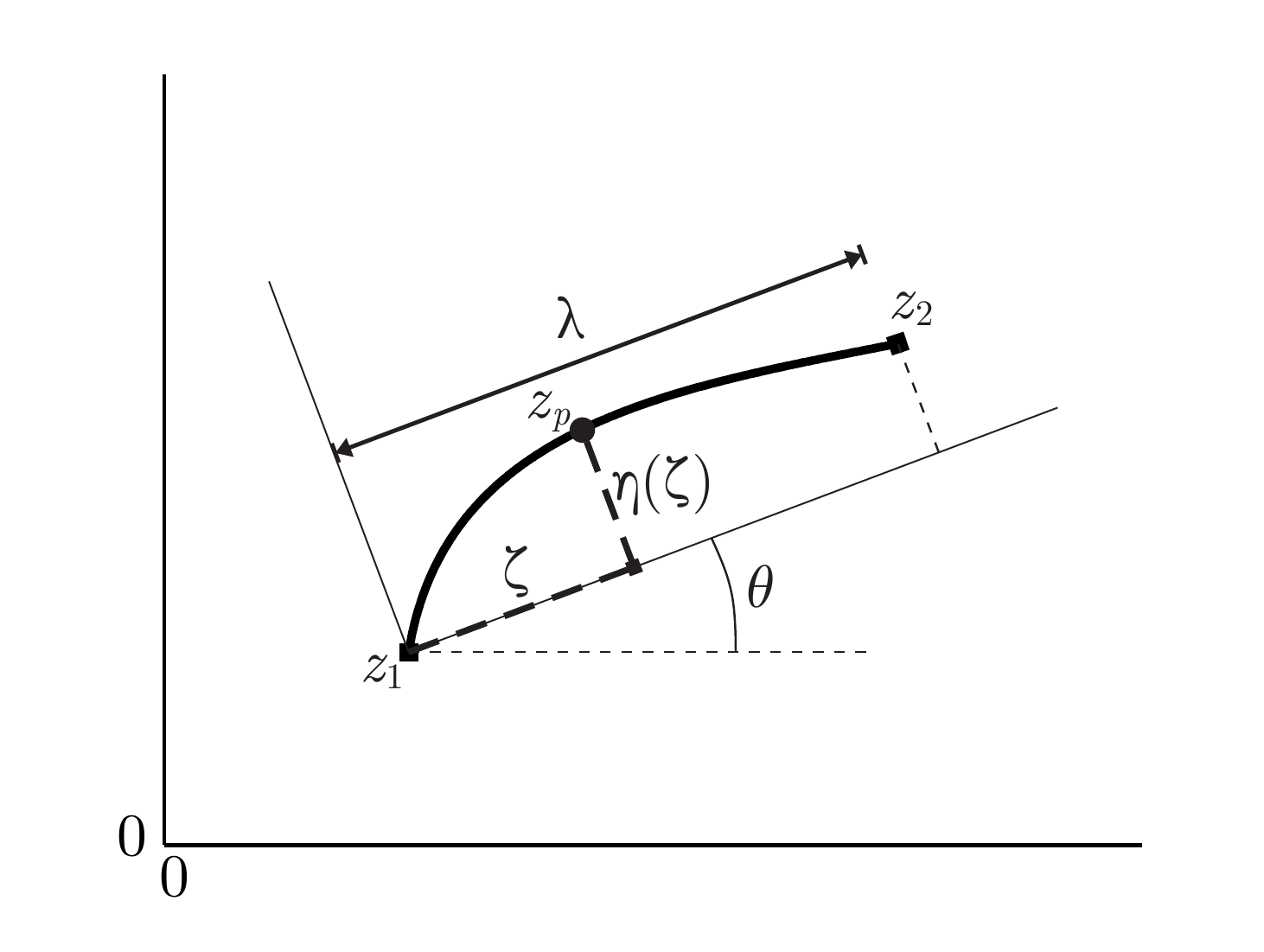}
  \caption{Schematics of a high order panel, showing both the global
    reference system, and the local reference system used in the
    derivations.}
  \label{fig:panelIllustration}
\end{figure}

Assume that a part of the boundary of an object extends between the
points $z_{1}$ and $z_{2}$. To generate a high order panel that models
this boundary, we need to select a baseline along which the panel is
parameterized. Since only the direction of this reference line will be
of importance, we let it start at $z_{1}$ with some angle $\theta$
pointing in the general direction of $z_2$, see
Figure~\ref{fig:panelIllustration}. The natural choice is to choose
the baseline between $z_{1}$ and $z_{2}$, but to allow for the panel
to be split (as relied upon in the FMM), it is necessary to allow that
$z_2$ is not on the baseline.

Following the notation of
\cite{Ramachandran03}, if we are interested in the flow velocity at
position $z$, we apply the transformation
\begin{equation}
  \label{eq:zprime}
  z' = z'(z) = \left(z-z_{1}\right)e^{-i\theta},
\end{equation}
which will make the reference line parallel to the real axis and the
panel will extend between 0 and $\lambda$, where
\begin{equation}
  \label{eq:lambda}
  \lambda=\mbox{Re}\left\{ \left(z_{2}-z_{1}\right)e^{-i\theta}\right\}.
\end{equation}
Using this reference system we can write a position $z_p$ on the panel as
\begin{equation}
  \label{eq:zp}
  z_{p} = z_p(\zeta) = 
  \zeta+i\eta\left(\zeta\right)=\zeta+i\sum_{k=0}^{M}a_{k}\zeta^{k},
\end{equation}
for a panel shape defined by a polynomial of degree $M$. To ensure
that the shape of the panel remains numerically reasonable, it will be
assumed that the coefficients $a_{k}$ are real.

The panel strength distribution is similarly defined by an $N$th
degree polynomial,
\begin{equation}
  \label{eq:gamma}
  \gamma = \gamma(\zeta) = \sum_{j=0}^{N}b_{j}\zeta^{j}.
\end{equation}
The coefficients $b_{j}$ are complex numbers, where the real part
represents the source strength and the imaginary part the vortex
strength according to \eqref{eq:vortexvelocity}.

\subsection{Evaluation of velocities}
\label{subsec:velocity}
The velocity contribution from the panel can be found by
integration,
\begin{align}
  \nonumber
  \overline{V\left(z\right)} &= \frac{e^{-i\theta}}{2\pi}\int\limits_0^\lambda
  \frac{\gamma(\zeta)}{\left(z'-z_{p}(\zeta)\right)} \,    
  d\zeta = \frac{e^{-i\theta}}{2\pi}\int\limits_0^\lambda
  \frac{\sum_{j=0}^{N}b_{j}\zeta^{j}}
       {\left(z'-\zeta-i\sum_{k=0}^{M}a_{k}\zeta^{k}\right)}
       \, d\zeta \\
  \label{eq:velc}
  &=\frac{e^{-i\theta}}{2\pi}\frac{-i}{a_{M}}\int\limits
  _{0}^{\lambda}\frac{\sum_{j=0}^{N}b_{j}\zeta^{j}}
        {\sum_{k=0}^{M}c_{k}\zeta^{k}} \, d\zeta,
\end{align}
where $z'$ is defined in \eqref{eq:zprime} and
\begin{align}
  c_{0} & = \frac{a_{0}+iz'}{a_{M}}, \quad c_{1} = \frac{a_{1}-i}{a_{M}},
          \quad c_{k} = \frac{a_{k}}{a_{M}}.
\end{align}
Hence the task can been reduced to an integration of the quotient of
two complex polynomials. First we ensure that the order of the
nominator is lower than the denominator (polynomial division),
\begin{equation}
  \label{eq:polydiv}
  \frac{\sum_{j=0}^{N}b_{j}\zeta^{j}}{\sum_{k=0}^{M}c_{k}\zeta^{k}} =
  \sum_{m=0}^{N-M}d_{m}\zeta^{m}+\frac{\sum_{j=0}^{M-1}h_{j}\zeta^{j}}{\sum_{k=0}^{M}c_{k}\zeta^{k}}
\end{equation}
Now, since $c_{M}=1$, we can factorize the denominator,
\begin{equation}
\sum_{k=0}^{M}c_{k}\zeta^{k}=\prod_{k=1}^{M}\left(\zeta-x_{k}\right).\label{eq:factorization}
\end{equation}
This allows a partial fraction decomposition, reducing the order of
the denominator,
\begin{align}
  \nonumber
  \frac{\sum_{j=0}^{N}b_{j}\zeta^{j}}{\sum_{k=0}^{M}c_{k}\zeta^{k}} &=
  \sum_{m=0}^{N-M}d_{m}\zeta^{m}+\frac{\sum_{j=0}^{M-1}h_{j}\zeta^{j}}{\prod_{k=1}^{M}\left(\zeta-x_{k}\right)} \\
  \label{eq:partialfrac}
  &= \sum_{m=0}^{N-M}d_{m}\zeta^{m}+\sum_{k=1}^{M}\frac{A_{k}}{\left(\zeta-x_{k}\right)},
\end{align}
where the constants $A_{k}$ are given by
\begin{equation}
  \label{eq:partialfraccoeff}
  A_{k}=\frac{\sum_{m=0}^{M-1}h_{m}x_{k}^{m}}{\prod_{n\ne k}\left(x_{k}-x_{n}\right)}.
\end{equation}
The final expression for the flow velocity from one panel thus becomes
\begin{equation}
  \label{eq:generalsolution}
  \overline{V\left(z\right)} =
  \frac{e^{-i\theta}}{2\pi}\frac{-i}{a_{M}}
  \left(\sum_{m=0}^{N-M}\frac{d_{m}}{m+1}\lambda^{m+1}+
  \sum_{k=1}^{M}A_{k}
  \log\left(\frac{x_{k}-\lambda}{x_{k}}\right)\right).
\end{equation}
This is the general expression for evaluating the flow velocity from
panels with arbitrary order polynomial shapes and strength
distributions. So far, the only numerical step is to find the roots of
the polynomial in the denominator, a well-studied problem for which
efficient algorithms exist \cite{Jenkins1970,Press:1992:NRC:148286}.

For numerical issues to consider when using this expression,
see~\ref{app:numerical}.

\subsubsection{Outgoing expansion}

Although \eqref{eq:generalsolution} provides a solution for the
velocity from an arbitrary polynomial panel, the expression is time
consuming to calculate and can also give rise to cancellation errors
when the distance $z'$ in \eqref{eq:zprime} is much larger than the
panel length $\lambda$. For a distance sufficiently far away from the
panel, a faster and more accurate way is to evaluate the velocity
through a series expansion. Let $z_{0}$ be a suitably chosen expansion
point for the panel. If we denote the radius of the smallest circle
that encloses the panel by $R$, then if $\left|z-z_{0}\right| > R$,
the evaluation can be carried out through the expansion
\begin{align}
   &\overline{V\left(z\right)} = \frac{1}{2\pi}\sum_{n=1}^{\infty}\frac{f_{n}}{\left(z-z_{0}\right)^{n}}\label{eq:Vseries},
\intertext{with}
\label{eq:fnseries}
f_{n}&=\int\limits _{0}^{\lambda}\sum_{j=0}^{N}b_{j}\zeta^{j}\left(\left(\zeta+i\sum_{k=0}^{M}a_{k}\zeta^{k}\right)e^{i\theta}-z_{0}+z_{1}\right)^{n-1}\,d\zeta.
\end{align}
Hence the coefficients $f_{n}$ can be evaluated through
straightforward integration of a polynomial. The number of
coefficients that has to be included in the series will depend on the
distance $\left|z-z_{0}\right|/R$ and on the desired precision.  Note
that this expansion is on the same form as the expansions used in the
fast multipole method. Hence, the panels can easily be included in the
multipole calculations by shifting these expansions to the center of
the FMM box containing the panel. Since the expansion is the same as
used in the FMM, the same error estimates can be used to determine the
required number of coefficients, see \cite{Greengard87} for details.
The current implementation uses these series expansions if $R < 1.8
\left|z-z_0\right|$, thus ensuring a relatively fast convergence of
the series.

\subsection{Solving the no-penetration boundary condition}
\label{subsec:bc}

The most common way of implementing the no-penetration boundary
condition is to define one control point at the center of each panel
and enforce that the normal flow velocity at this point is zero. For
some applications it is desirable that the flow does not leak out of
the domain, and it can be more suitable to satisfy the boundary
condition that the net flow through the entire panel should be zero,
\begin{align}
\int\limits _{Panel}\vec{V}\cdot\hat{n} \, dS=0.
\end{align}
This condition also makes the solution less sensitive to vortices
located close to a control point. 

\subsubsection{Contribution from a point source}
The contribution of one single
vortex on a panel can be determined as follows. Start by applying
Gauss' law,
\begin{align}
\oint\limits _{S}\vec{V}\cdot\hat{n}dS=\iint\limits _{A}\nabla\cdot\vec{V} \, dA.
\end{align}
As a closed contour $S$ is required, we close the contour by simply
inserting a straight line between the end points. That is,
\begin{align}
\oint\limits _{S}\vec{V}\cdot\hat{n} \, dS=\int\limits _{Panel}\vec{V}\cdot\hat{n} \, dS+\int\limits _{Line}\vec{V}\cdot\hat{n} \, dS,
\end{align}
which will give us the flow through the panel as
\begin{equation}
\int\limits _{Panel}\vec{V}\cdot\hat{n}\,dS=-\int\limits _{Line}\vec{V}\cdot\hat{n} \, dS+\sum_{q\in A} q,\label{eq:panelgauss}
\end{equation}
where $q$ represents the strength of the vortex and where $A$ is the
area enclosed by $S$. This shows that the integral only has to be
calculated along the line, and then a correction is applied whenever
the vortex is situated within $A$. In turn, the line integral is given
by
\begin{equation}
  \int\limits _{Line}\vec{V}\cdot\hat{n}\,dS=\mbox{Re}\left(\frac{q}{2\pi}\log\left(\frac{z'-\lambda}{z'}\right)\right)=\mbox{Re}\left(\frac{q}{2\pi}\log\left(\frac{z-z_{2}}{z-z_{1}}\right)\right).\label{eq:panellinevel}
\end{equation}
Note that the imaginary part of the expression is the flow parallel to
the panel.

If so-called smoothing kernels of the vortices are used
\cite{Cottet08}, and if the kernel overlaps with the panel, numerical
integration may be the best option.

\subsubsection{Evaluation from local- and far series expansions}

In the fast multipole method, each multipole box will be associated
with a local field expansion of the form
\begin{align}
\overline{V\left(z\right)}=\sum_{n=0}^{P}p_{n}\left(z-z_{0}\right)^{n},
\end{align}
where $z_{0}$ is the center of the expansion. This expansion contains
the contribution from all vortices far away, and it is hence much
faster to evaluate the net panel flow through this expansion rather
than from \eqref{eq:panellinevel}. If we carry out the integral over
the line, we find
\begin{align}
\int\limits _{Line}\vec{V}\cdot\hat{n}dS=\sum_{n=0}^{P}p_{n}\frac{\left(z_{2}-z_{0}\right)^{n+1}-\left(z_{1}-z_{0}\right)^{n+1}}{n+1}.
\end{align}

When evaluating the contribution from another panel far away, yet
close enough to be in the near-field of the fast multipole method, it
is preferable to carry out the evaluation from the far-field
expansion around some point $z_{0}$. The general form of these
expansions is given in \eqref{eq:Vseries}, here truncated to $P$
coefficients.  We find
\begin{align}
  \int\limits _{Line}\vec{V}\cdot\hat{n}\,dS = \Biggl[f_{1}\log\left|z-z_{0}\right|-\sum_{n=2}^{P}\frac{f_{n}}{\left(n-1\right)\left(z-z_{0}\right)^{n-1}}\Biggr]_{z = z_{1}}^{z = z_{2}}.
\end{align}
This expression can be used when \eqref{eq:Vseries} is valid over the
line of integration.

\subsubsection{Panel to panel interaction}

To use the flow through the panel when generating the dense BEM
interaction matrix, we need to integrate the flow from one panel over
another panel. This expression can be obtained by integrating
\eqref{eq:panellinevel} over the source panel. By rewriting the panel
the same way as in \eqref{eq:velc}, we find
\begin{align}
  \nonumber
  \int\limits _{0}^{\lambda}\frac{1}{2\pi}\gamma\log\left(\frac{z-z_{2}}{z-z_{1}}\right) &= 
  \frac{1}{2\pi}\frac{-i}{a_{M}}\int\limits _{0}^{\lambda}\sum_{j=0}^{N}b_{j}\zeta^{j}\log\left(\sum_{k=0}^{M}c_{2,k}\zeta^{k}\right)d\zeta\\
  \label{eq:paneptopanelbase}
  -\frac{1}{2\pi}&\frac{-i}{a_{M}}\int\limits _{0}^{\lambda}\sum_{j=0}^{N}b_{j}\zeta^{j}\log\left(\sum_{k=0}^{M}c_{1,k}\zeta^{k}\right)d\zeta+N_{1}i\int\limits _{0}^{\lambda}q.
\end{align}
The constant $N_{1}$ has been included to compensate for the possible
changes in the branch of the logarithm. Equation
\eqref{eq:paneptopanelbase} shows that the solution is reduced to the
solution of two logarithmic potentials. A brief derivation of the
solution to the logarithmic potential is included in~\ref{app:logpot}.

Using the expression for the logarithmic potential and by assuming
that the panels do not intersect each other, the panel to panel
interaction can be written as
\begin{align}
\nonumber
\frac{1}{2\pi}\int\limits _{0}^{\lambda}\gamma\log\left(\frac{z-z_{2}}{z-z_{1}}\right) = & \frac{1}{2\pi}\left\{ \sum_{j=0}^{N}\frac{b_{j}}{j+1}\left(\lambda^{j+1}\right)\log\left(\frac{\sum_{k=0}^{M}c_{2,k}\lambda^{k}}{\sum_{k=0}^{M}c_{1,k}\lambda^{k}}\right)\right.\\
 & \left.+\sum_{k=1}^{M}\left(\sum_{m=0}^{N}\left(s_{m}\left(x_{2,k}\right)-s_{m}\left(x_{1,k}\right)\right)\lambda^{m}\right)\right\}
\end{align}
where
\begin{align*}
s_{0}\left(x_{k}\right) &:= -\sum_{j=0}^{N}\frac{b_{j}}{j+1}\left(x_{k}^{j+1}\right)\log\left(\frac{x_{k}-\lambda}{x_{k}}\right),\\
s_{m}\left(x_{k}\right) &:= \frac{1}{m}\sum_{j=m-1}^{N+1}b_{j}\frac{1}{\left(j+1\right)}x_{k}^{j+1-m}.
\end{align*}
Here, $x_{1,k}$ and $x_{2,k}$ are the roots to the polynomials defined
by coefficients $c_{1,k}$ and $c_{2,k}$ respectively.

The above expression experiences a problem when the end point $z_{2}$
of the source panel coincide with any of the end points of the
evaluation panel. If this is the case, the easiest solution is to
generate a new panel which is parameterized in the opposite direction,
hence moving the root from $\lambda$ to 0.

A correction has to be added whenever the panel end up within the area
that is enclosed by the line between the evaluation panel end points
and the evaluation panel itself. In this case, a constant with the
magnitude
\begin{equation}
 C=\pm i\int\limits _{0}^{\lambda}q
\end{equation}
 should be added (see the last term in \eqref{eq:paneptopanelbase}),
 where the sign is chosen depending on the direction of the integral
 over the line and in relation to the direction in
 \eqref{eq:panelgauss}.

\subsection{Continuity of the source strength}
\label{subsec:smooth}

As stated in \S\ref{sec:Introduction}, to generate a smooth velocity
field, the shape of the panels should have a continuous derivative and
the strength should be continuous. In the current definition of the
panels, however, the source strength is distributed along the
reference line rather than along the tangential direction of the
panel. It follows that the strength of the panel (in the tangential
direction of the panel) will depend on the slope of the panel. To
ensure continuity when determining the panel strength, it is
convenient to define the panel by its strength (and derivatives of the
strength for higher order continuity) in the end points. Hence, we
want to give the strength in the end points along the tangential
direction of the panel when defining the panel.  One possible approach
is to modify \eqref{eq:gamma} to read
\begin{equation}
  \label{eq:gammaactual}
  \gamma = \gamma(\zeta) = 
  \sum_{j=0}^{N}b_{j}\zeta^{j}\sqrt{1+\left(\frac{d\eta}{d\zeta}\right)^{2}}.
\end{equation}
Assuming that the slope of the panel is small, it is possible to
expand \eqref{eq:gammaactual} into a series, effectively a polynomial,
where \eqref{eq:generalsolution} applies anew. Any panels with large
slope can be split into smaller panels until the series expansion is
convergent. The disadvantages of this approach is that it generates
additional panels and that the resulting polynomial in the nominator
of \eqref{eq:velc} is of high degree.  As an alternative approach, one
can choose to only apply a correction factor to the given input values
for the strength at the panel end points.  Assume that the desired
strength in the tangential direction of the panel is given by
$\gamma_{input}$, then one can apply the correction at this point
\begin{equation}
\gamma=\gamma_{input}\sqrt{1+\left(\frac{d\eta}{d\zeta}\right)^{2}},\label{eq:gammacorr}
\end{equation}
and then use \eqref{eq:gamma} without any modification. This will give
a different solution, than using \eqref{eq:gammaactual}, but
continuity will still be fulfilled. The corrections for higher order
continuity are obtained by differentiating \eqref{eq:gammacorr} with
respect to $\zeta$.


\section{Numerical experiments}
\label{sec:evaluation}

This section will illustrate the performance of the higher order
panels for two different cases. In all cases, the panels are specified
by giving the position, the derivative, the second order derivative
and so on at each end point, where the number of derivatives are
chosen depending on the order of the panel. We have tested panel
shapes of all odd orders up to the 7th order. The panel strength is
also chosen by giving the value, derivative, etc.\ at each end point,
hence only odd orders are included here as well. This also means that
panels of order 3 will have continuous first derivatives, panels of
order 5 will additionally have continuous second derivatives
etc. Consequently, 1st order panels imply 1 unknown per panel, 3rd
order panels 2 unknowns, 5th order 3 unknowns and so on. All
simulations will apply the correction to the panel strength at the end
points according to \eqref{eq:gammacorr} and its derivatives.

If more unknowns are added, there will also be a need for satisfying
the zero flow condition on more panels. To generate more panels, the
source panels are split into several evaluation panels while
maintaining the same shape. Also, due to issues with BEM matrices that
are close to singular, in the general case, it was concluded that
better solutions are obtained if the source panels are split into one
evaluation panel too much, and then the zero flow conditions is solved
in a least square sense. All solutions in this section will apply this
technique.

\subsection{Flow around a deformed circle}

\label{sub:unitcircle}
\begin{figure}
  \begin{center}
    \includegraphics[width=0.7\columnwidth]{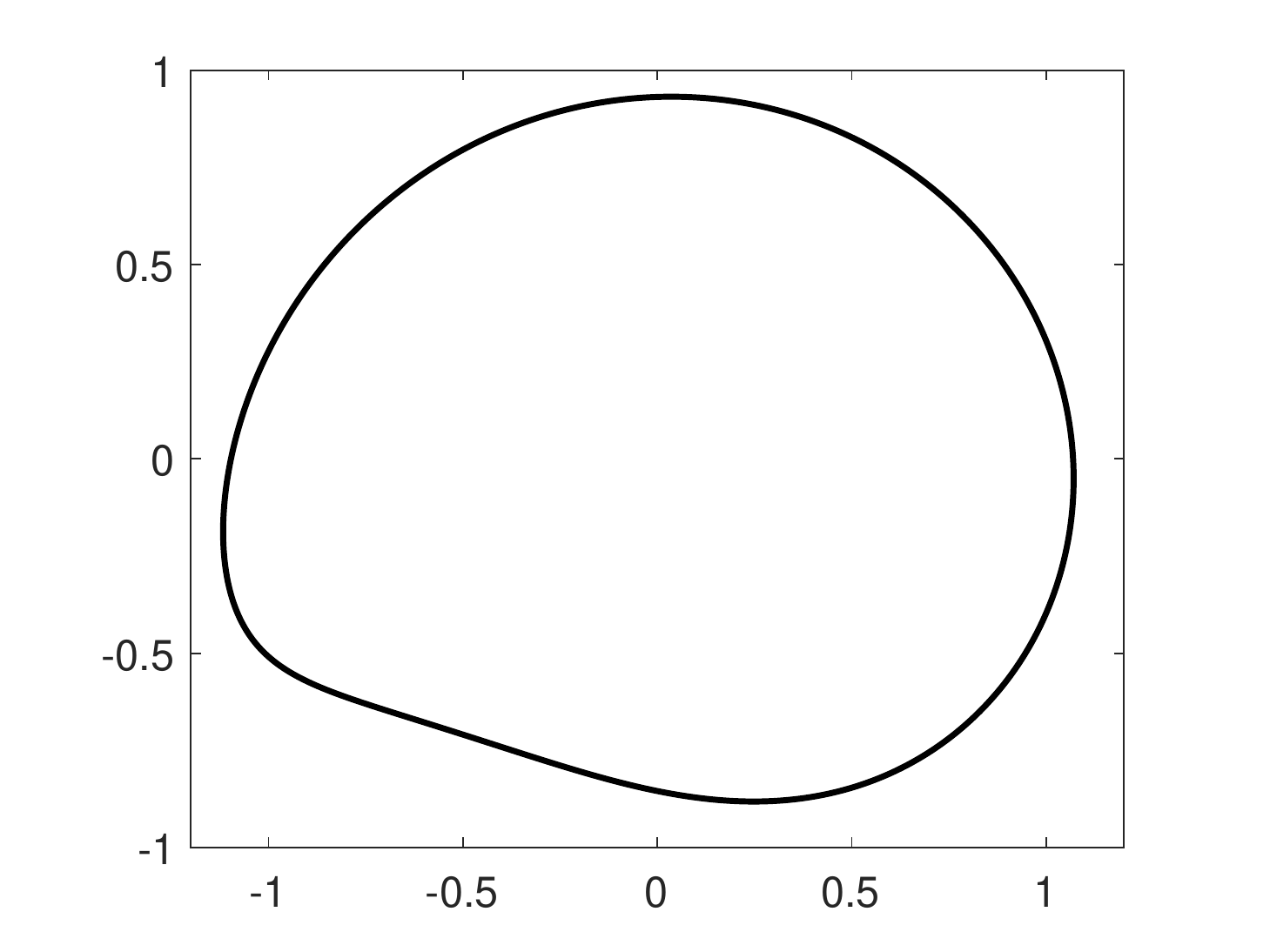}
    \caption{\label{fig:deformed_circle}The deformed circle used in
      the tests.}
  \end{center}
\end{figure}
\begin{figure}
  \includegraphics[width=1\columnwidth]{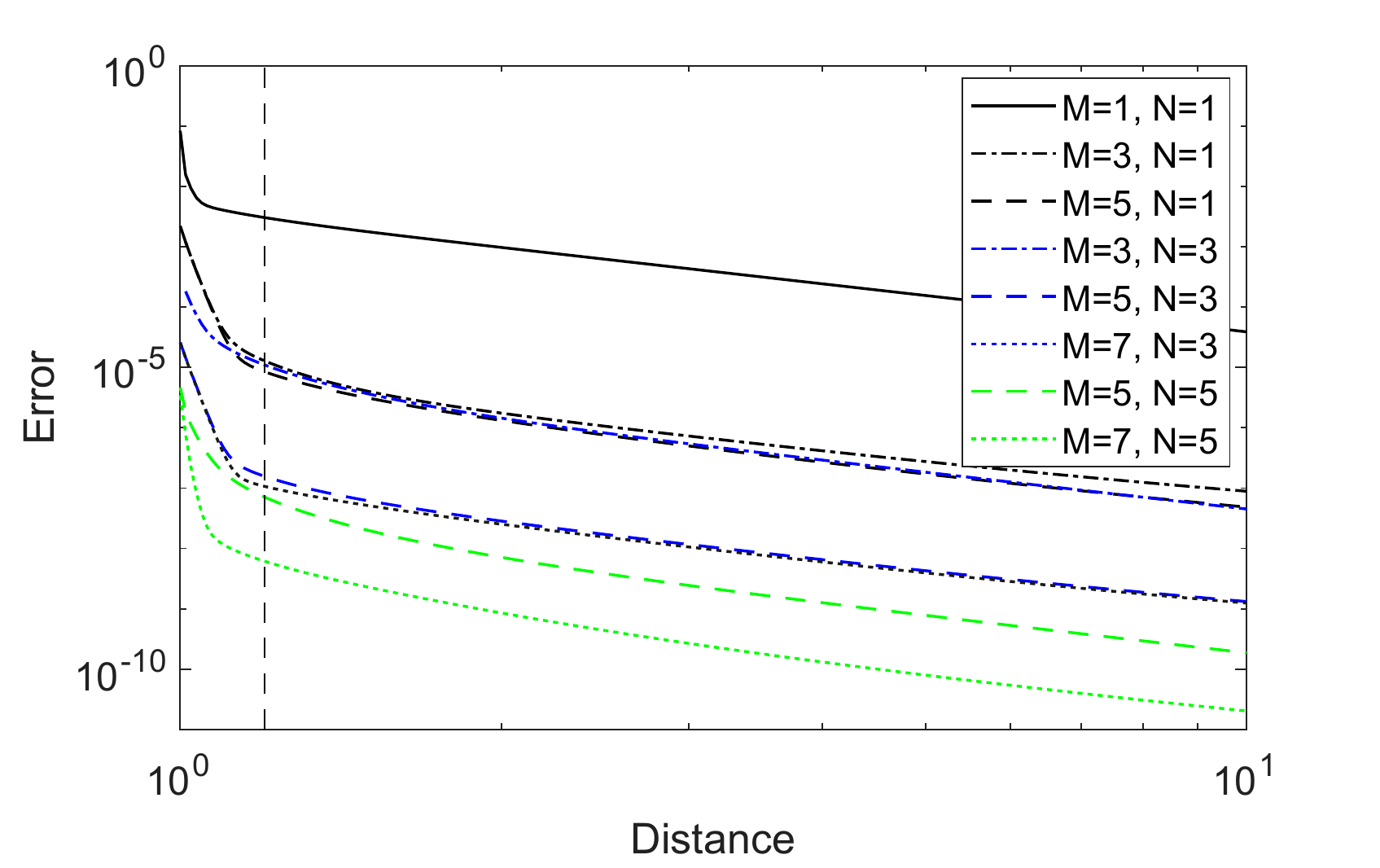}
  \caption{\label{fig:circledistanceconvergence}Error as function of
    distance from the deformed circle for 50 panels. Note that the
    surface of the deformed circle is located at distance 1. In the
    legend, $M$ is the order of the shape polynomial and $N$ is the
    order of the strength polynomial. The vertical dashed line shows
    the distance 1.2, which is used in
    Figure~\ref{fig:circleNpanelconvergence}.}
\end{figure}
\begin{figure}
  \includegraphics[width=1\columnwidth]{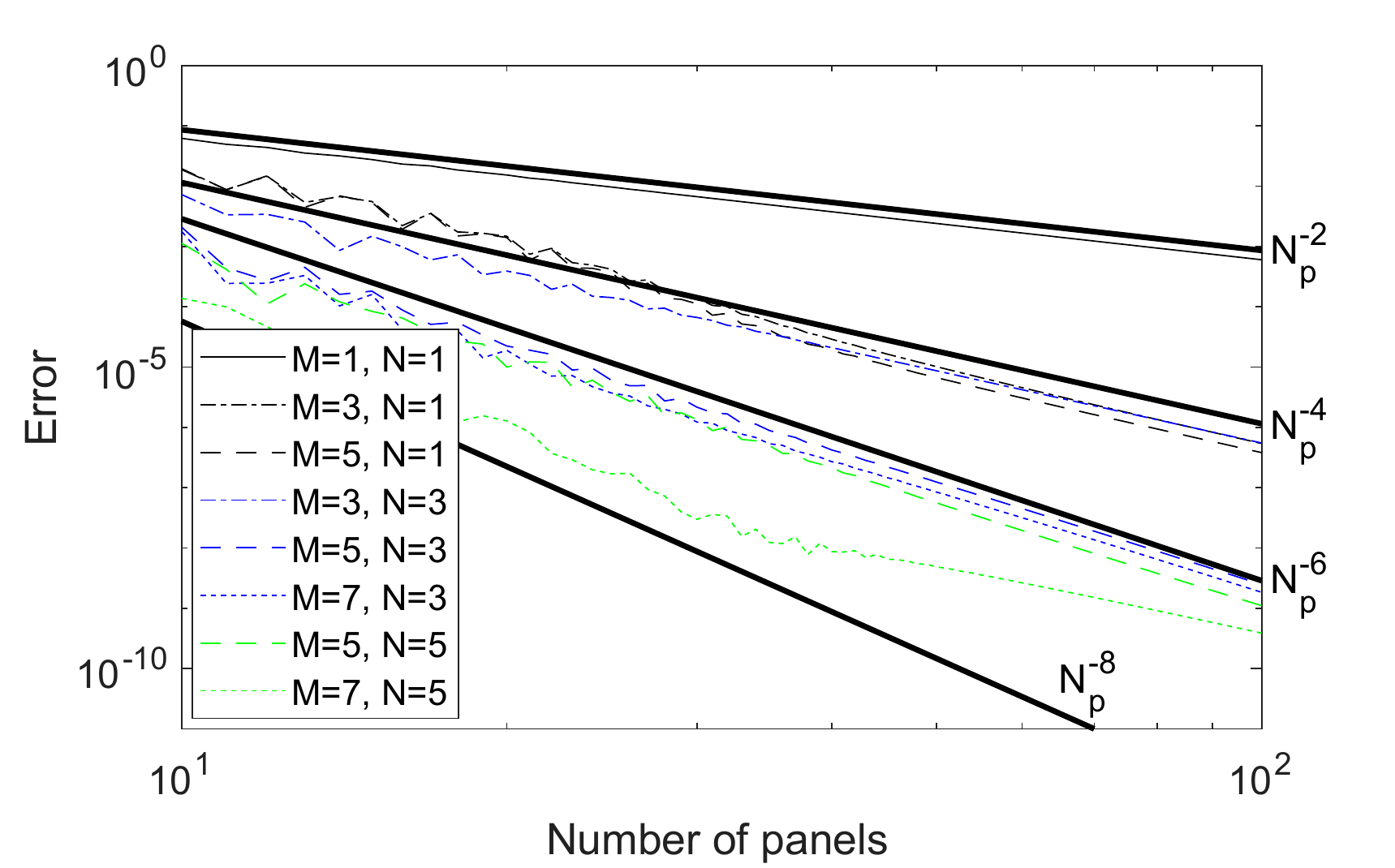}
  \caption{\label{fig:circleNpanelconvergence}Convergence as a
    function of the number of panels. The bold black lines represent
    different orders of convergence. }
\end{figure}
For the first test case, a smooth shape in the form of a deformed
circle is chosen. The deformation is accomplished by using the
conformal mapping
\begin{equation}
  z=s+\frac{0.1}{s+0.3+0.4i}
\end{equation}
of a unit circle, and the resulting shape is illustrated in
Figure~\ref{fig:deformed_circle}. The deformation is added to the
circle to avoid testing the panels on the symmetric circle, which
artificially can give higher convergence than the general case.

The flow velocity when a conformal mapping is applied to a unit circle
(without circulation) has the well-known analytic solution
\begin{equation}
\overline{V_{ref}\left(z\right)}=\frac{V_{\infty}\left(1-\frac{1}{s^{2}}\right)}{\frac{dz}{ds}}.
\end{equation}
The measurement of error that will be used here is 
\begin{equation}
  E=\frac{1}{\sqrt{L}V_{\infty}}\sqrt{\int\limits _{S}\left(V_{ref}\left(z\right)-V\left(z\right)\right)^2dz},
  \label{eq:circleerror}
\end{equation}
where $V$ is the velocity calculated from the panels, $S$ is a curve
around the deformed circle at a fixed distance from the object and $L$
is the length of the curve $S$. This expression is evaluated through
numerical integration with Gauss-Kronrod quadrature in
Matlab~\cite{Shampine2008131}.

For this case, panels with pure source strength are chosen (i.e.~no
vorticity) and the additional criterion that the total source strength
of all panels should be zero is also enforced. All panels are
distributed evenly around the unit circle in the $s$-plane.

As a first evaluation, the error is given in
Figure~\ref{fig:circledistanceconvergence} as a function of the
distance from the deformed circle when using 50 panels. Some distance
away from the circle (in the far-field of the panels), the error
decreases with approximately the same factor for all panels, but the
higher order panels give significantly better results for the same
amount of panels. Very close to the circle (starting when the distance
from the circle is about the same as the length of the panel), there
is a significant change in the slope of the error curves and the error
grows as the distance approaches the object.  One general trend is
that the order of the strength of the panels appears to be more
important close to the surface of the object.

For linear strength ($N=1$), compared to the standard approach of
using linear panels ($M=1$), there are significant improvements for
increasing to cubic shape, but very modest improvements for the
quintic shape. For cubic strength, a quintic shape is required to
obtain significant improvements compared to the linear strength (while
septic shape gives modest improvements) and a quintic shape shows some
improvements for both the quintic and the septic shape.  These trends
are also visible in Figure~\ref{fig:circleNpanelconvergence}, where
the convergence in the far-field with respect to the number of panels
is shown. Here, it can be seen that panels with linear shape have 2nd
order convergence (regardless of the order of the strength), panels
with cubic shape have 4th order convergence if strength is linear or
higher, quintic panels with have 6th order convergence if strength is
cubic or higher and septic panels initially have 8th order convergence
if strength is quintic or higher.  The trend change for the septic
panels are likely related to numerical errors in the solution of the
BEM equations. The maximum possible convergence rates in the far-field
is according to the derivation in~\ref{app:convergence} equal to
$\min(M+1,2N+2)$ under the condition that the optimal strength
coefficients are obtained when solving the BEM equations. However,
this assumes that the number of unknowns for the BEM equations is
equal to $N+1$. Here, higher order continuity is enforced, which
reduces the number of unknowns by a factor 2. Second, the BEM
equations are applied in the near-field, and solves the no-penetration
boundary condition, while the error is measured as the least square
error in the far-field in
Figure~\ref{fig:circleNpanelconvergence}. The trend seen here is that
the theoretical estimate applies to constant and linear strength,
while when increasing the strength to cubic or higher, the reduction
by a factor 2 applies (panels with constant strength are not shown in
Figure~\ref{fig:circleNpanelconvergence} but have quadratic
convergence with similar performance as $M=1$, $N=1$ for all tested
shapes).

For the deformed circle case, the conclusion is that significant
improvements can be gained from increasing the order of the panels. It
is also clear that having a good approximation of the shape is more
important in this case, than a good approximation of the strength.

\subsection{Flow around an airfoil}

\begin{figure}
  \includegraphics[width=1\columnwidth]{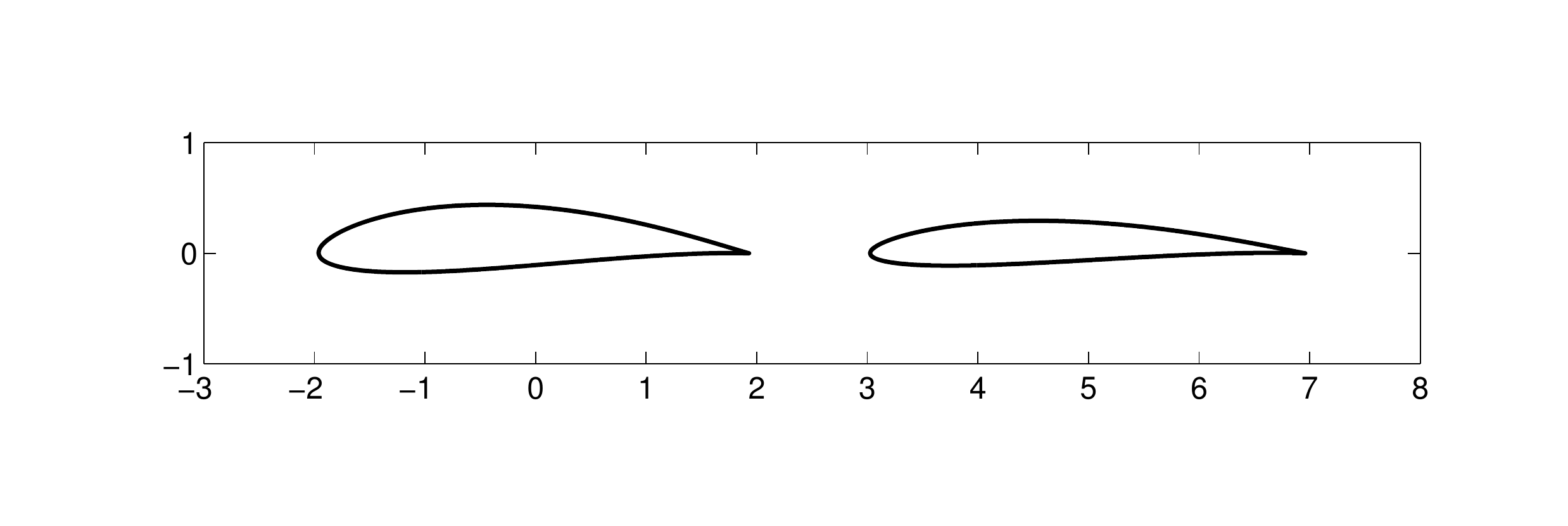}
  \caption{\label{fig:airfoils}The two tested airfoils. Airfoil to the
    left (A) is generated with $\mu=-0.09+0.09i$ and $n=1.93$, and
    airfoil to the right (B) is generated with $\mu=-0.06+0.06i$ and
    $n=1.95$ }
\end{figure}
In the second reference case, the flow around an airfoil is studied,
which is a common application for both panel methods and vortex
methods. Here, two aspects are important. The first is how well the
flow close to the boundary is approximated. This is important for
solutions of the boundary layer and therefore have many applications
\cite{Zanon13,Riziotis08}. The second aspect is how well the
circulation around the airfoil is calculated, which will determine the
calculated lift force. For this case, a Karman Trefftz airfoil is
chosen, which has an analytic solution through the use of the
conformal mapping
\begin{equation}
z=n\frac{\left(1+\frac{1}{s}\right)^{n}+\left(1-\frac{1}{s}\right)^{n}}{\left(1+\frac{1}{s}\right)^{n}-\left(1-\frac{1}{s}\right)^{n}},
\end{equation}
which maps a cylinder in the complex $s$-plane with center in $\mu$
and radius $R_{ref}=\left|1-\mu\right|$ to an airfoil. The velocity
for zero degrees pitch angle is obtained from
\begin{equation}
\overline{V_{ref}\left(z\right)}=\frac{V_{\infty}\left(1+\frac{i\Gamma_{ref}}{2\pi\left(s-\mu\right)}-\frac{R_{ref}^{2}}{\left(s-\mu\right)^{2}}\right)}{\frac{dz}{ds}},
\end{equation}
where 
\begin{equation}
\Gamma_{ref}=4\pi V_{\infty}\mbox{Im}\left(\mu\right)
\end{equation}
is the circulation. Here, the error in circulation is defined as
\begin{equation}
E_{circ}=\frac{\Gamma-\Gamma_{ref}}{\Gamma_{ref}}
\end{equation}
and the error in flow velocity is calculated according to
\eqref{eq:circleerror} by integrating along a line with with a fixed
distance from the surface (an arc is added around the trailing edge to
connect the upper and lower sides). The circulation around the airfoil
is calculated through the Kutta condition~\cite{Abbott59} by setting
the vortex strength to zero at the trailing
edge~\cite{drela2014flight}.

\begin{figure}
  \includegraphics[width=1\columnwidth]{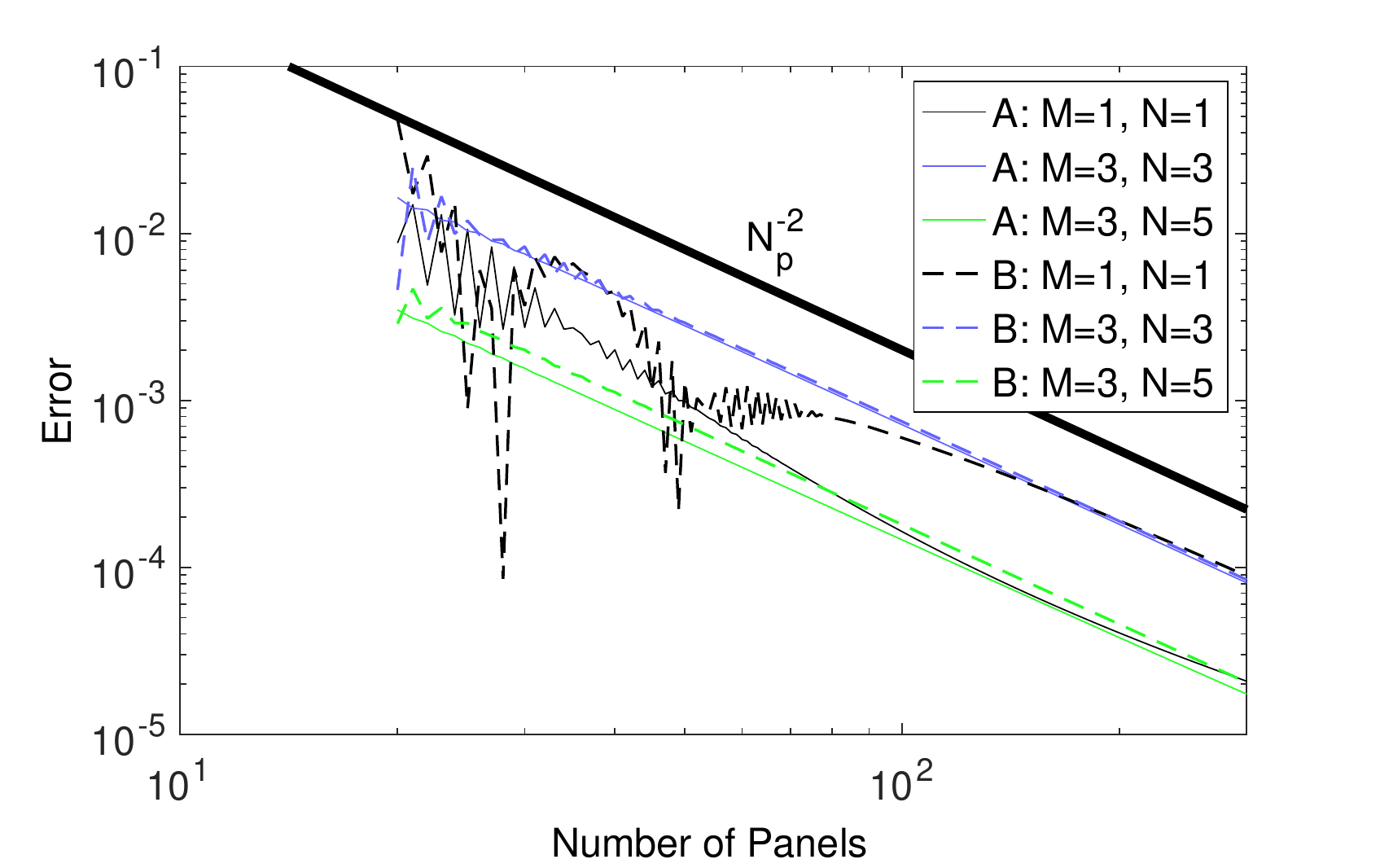}
  \caption{\label{fig:circulation}Error in circulation for airfoils A
    and B. Only panels up to cubic shape ($M=3$) are shown, as the
    higher order panels give almost identical result.}
\end{figure}

Here, two different profiles (A and B) are tested, see
Figure~\ref{fig:airfoils}. The panels are generated by evenly
separating the points in the $s$-plane with start and end panels on
the trailing edge. This will generate smaller panels close to the
trailing edge, which is the most difficult part to solve.  The BEM
matrix is then normalized to give all evaluation panels the same
weight before the least square fit (which makes the flow close to the
trailing edge more important). As the circulation is of interest, the
airfoils are modeled with pure vortex panels.

For numerical reasons, defining higher order panels through their
higher order derivatives does not work properly for quintic panels, as
the second order derivative grows close to the trailing edge. Hence
the order of shape for all points above 1.9 in the $x$-coordinate has
to be reduced to cubic. Note that this only affects a very small part
of the panels close to the trailing edge.

The calculated circulation for the two airfoils are shown in
Figure~\ref{fig:circulation}. This shows that the convergence is
similar for all panels, hence the high order convergence obtained with
the circle is not present here. This can probably be related to the
large gradients of the solution close to the trailing edge, which
causes numerical difficulties. The oscillations seen for the linear
panels are related to the modeling of the leading edge.  The high
accuracy values are obtained when the leading edge panel is
perpendicular to the incoming flow, while the low accuracy values are
obtained when the transition between two panels are at the leading
edge, giving a sharp edge there.  It should be noted that it is not
the same panel orders that give the best results in both cases. Pure
linear panels do for example give among the best results for airfoil
A, but not for airfoil B. Apparently, determining the best panel order
\textit{a priori} is not so easy.

\begin{figure}
  \includegraphics[width=1\columnwidth]{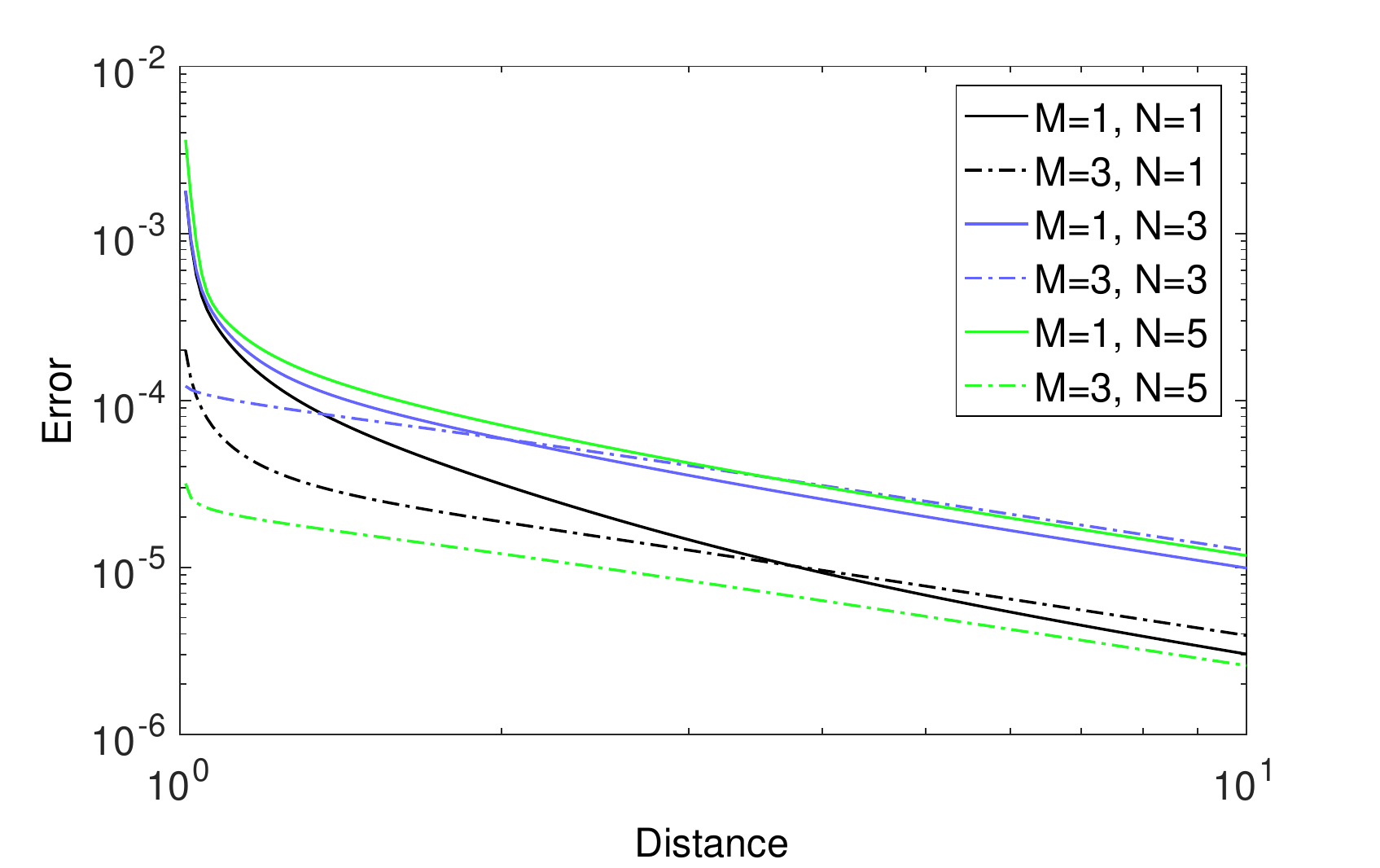}
  \caption{\label{fig:distanceA}Error in flow velocity as function of
    distance from the panel surface for 100 panels. The panel surface
    is located at distance 1.}
\end{figure}

\begin{figure}
  \includegraphics[width=1\columnwidth]{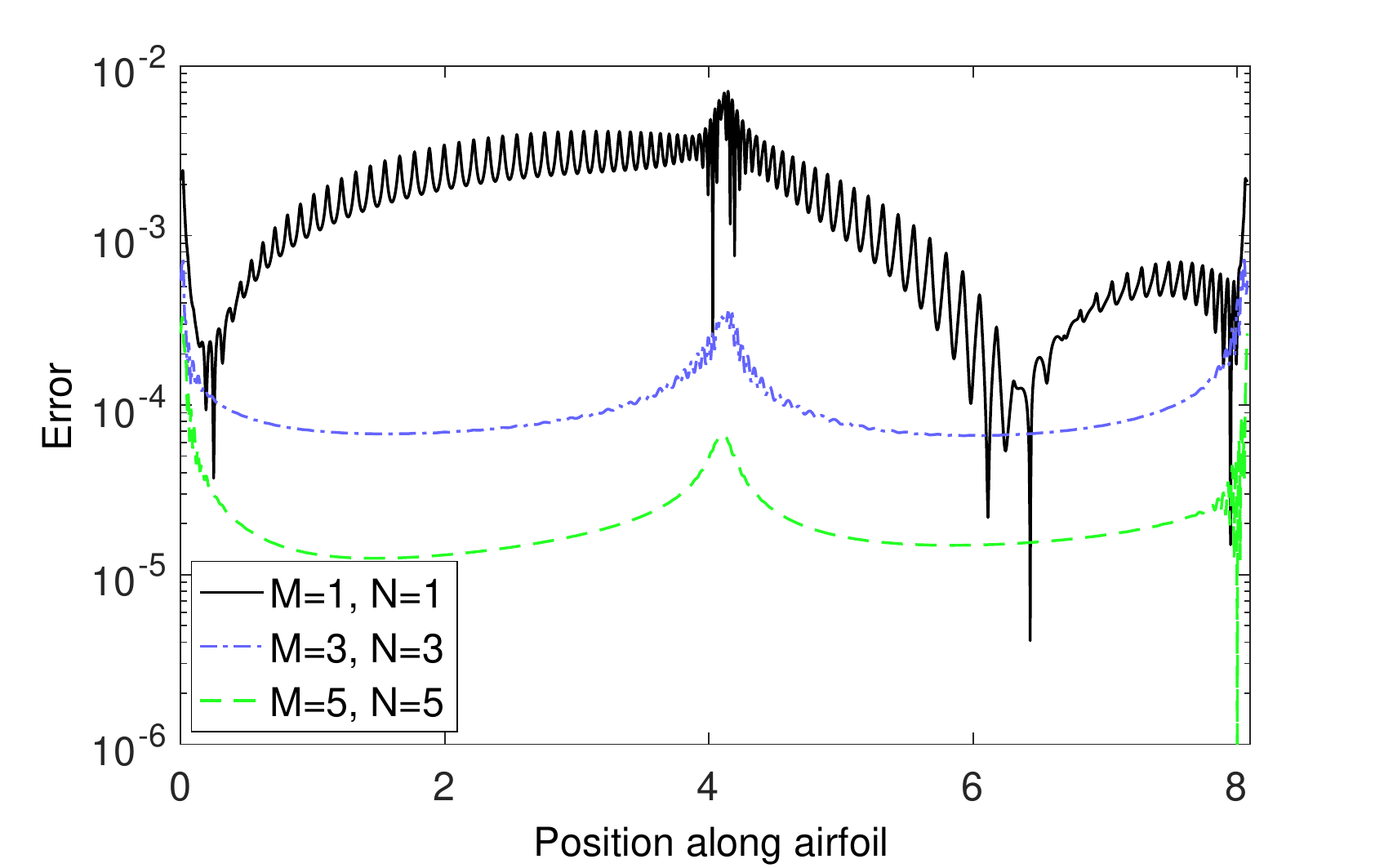}
  \caption{\label{fig:surfaceA}Error in flow velocity as a function of
    position along the airfoil surface at a distance of 0.01 from the
    surface for airfoil A with 100 panels. The position starts at the
    trailing edge and moves counter clockwise around the airfoil,
    ending at the trailing edge.}
\end{figure}

The convergence with respect to distance (Figure~\ref{fig:distanceA})
shows similar trends as for the circle. Some distance away, most panel
orders tend to have the same convergence. The exception is panels with
linear shape and strength, which show significant improvements with
respect to distance. This probably comes from the fact that the linear
panels obtained a good value for the circulation, and having the
correct total circulation becomes more important with distance. It is
also clearly seen that close to the panel surface, the higher order
panel shapes show significantly improved results. To further
illustrate the origin of this, the error as a function of distance
along the integration curve is shown in Figure~\ref{fig:surfaceA} at a
distance of 0.01 from the airfoil surface. Here, panels with linear
shape have large errors over the entire airfoil surface, while the
higher order panels show good values over most of the surface, but
with large errors close to the trailing edge (and somewhat larger
errors close to the leading edge).

\begin{figure}
  \includegraphics[width=1\columnwidth]{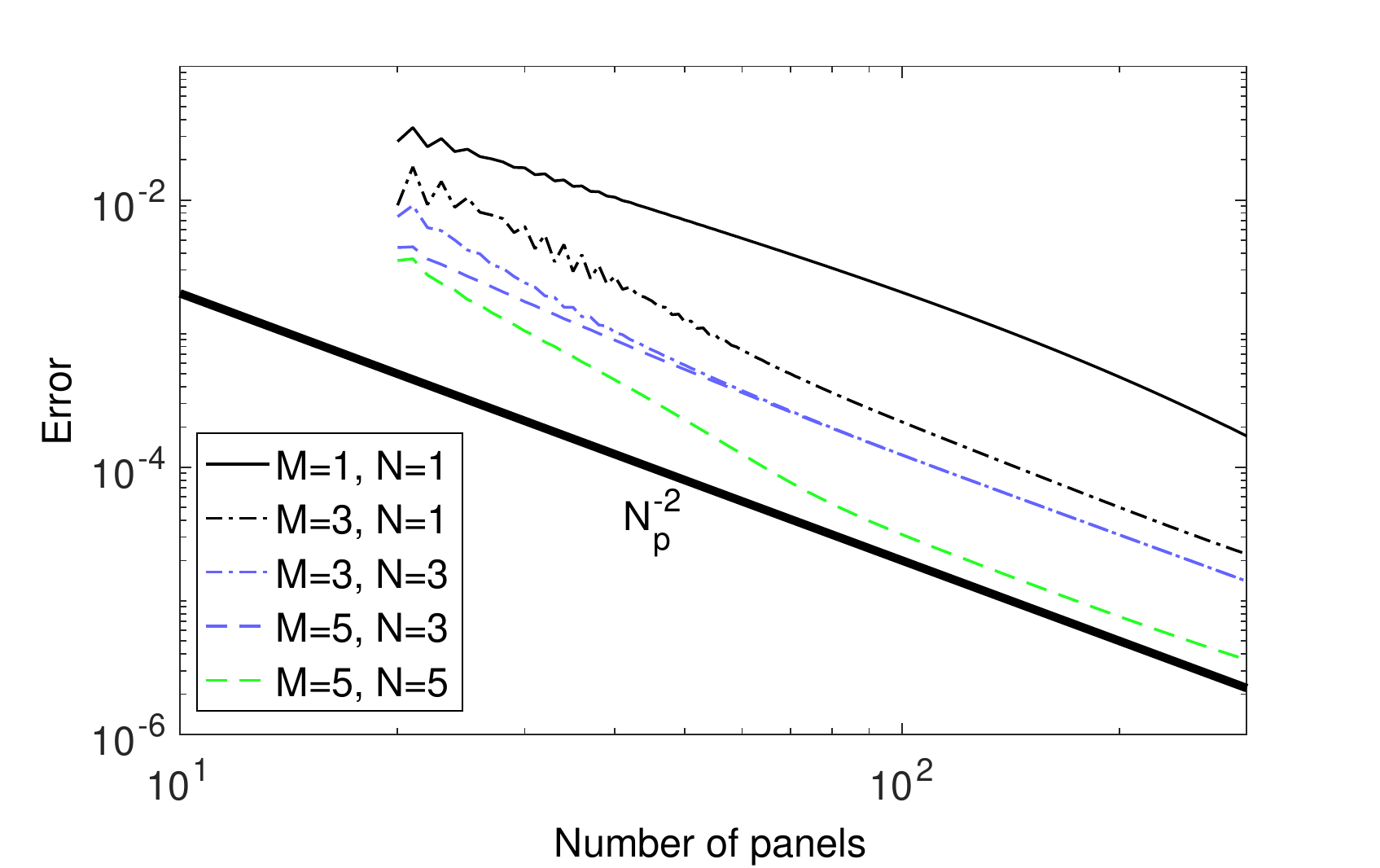}
  \caption{\label{fig:NconvergenceA}Error in flow velocity as a
    function of number of panels for airfoil A at a distance of 0.01
    from the surface.}
\end{figure}

Finally, the error in velocity at distance 0.01 as a function of the
number of panels is illustrated in
Figure~\ref{fig:NconvergenceA}. Similar to the error in circulation,
all different panel orders tend to have relatively similar convergence
behavior (compared to the large differences for the circle). Again,
the quintic panels with quintic strength show the best results, but
only quadratic convergence is obtained. Panels of septic order were
tested as well, but without any improvements to the results.

\subsubsection{Comparison of different panel implementations}

\begin{figure}
  \includegraphics[width=1\columnwidth]{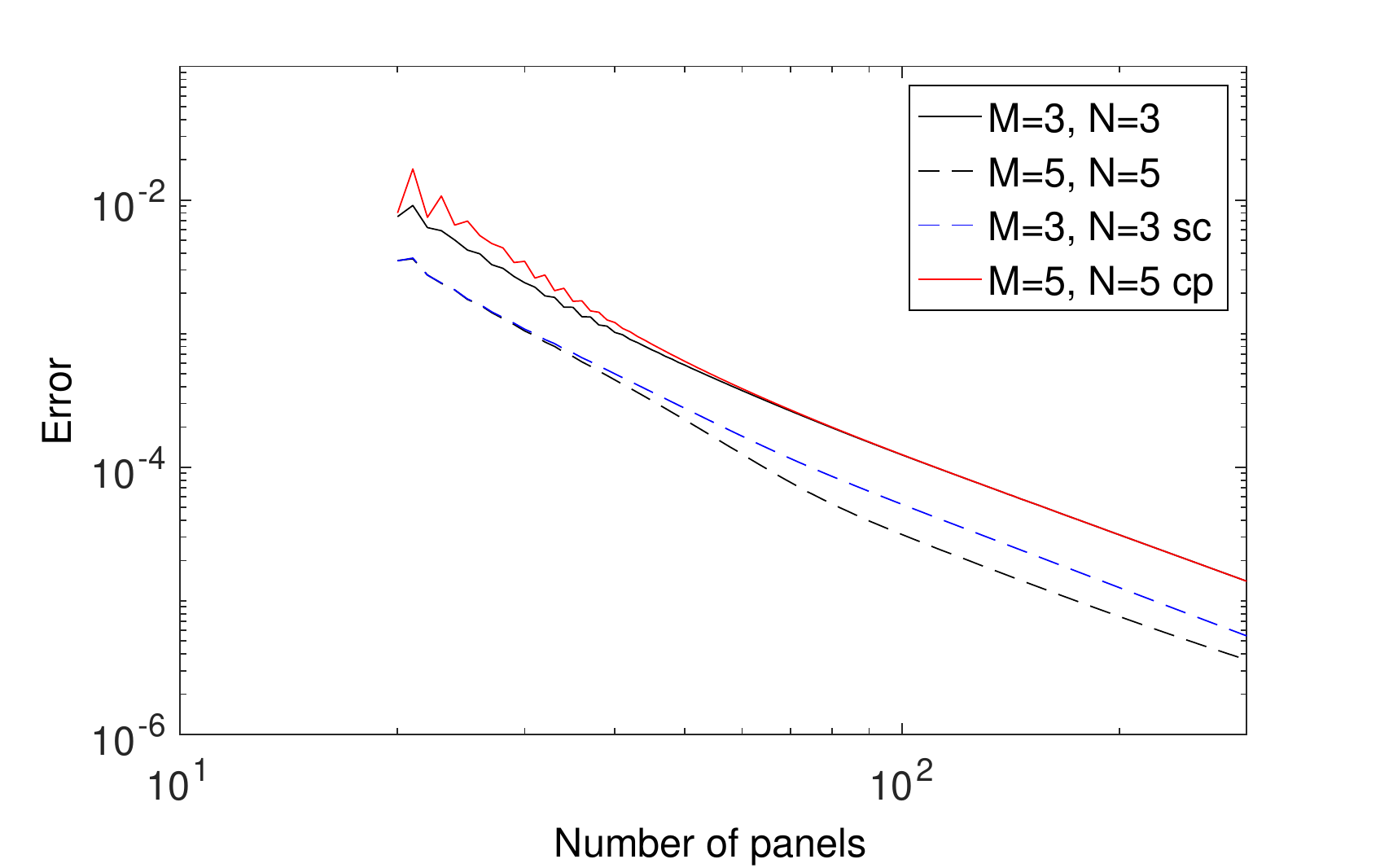}
  \caption{\label{fig:implementations}Comparison of the main method,
    with the use of control points in the center of each panel
    (instead of integrating the flow through the panel) marked with
    ``cp'' in the figure (blue line) and comparison of using corrected
    panel strength as in \eqref{eq:gammaactual}, marked as ``sc'' in
    the figure (red line), instead of only correcting the end points
    as in the main method.}
\end{figure}

One question is how big the difference is between solving for zero
flow through the panels, or solving for zero flow velocity at control
points in the center of the panels. To investigate this, a comparison
between the two have been added in Figure~\ref{fig:implementations}
for airfoil A (error in flow velocity at distance 0.01). A small
increase in performance is seen for quintic panels when using the flow
through the panel while for cubic panels, no difference was seen. This
shows that a control point approach is a decent approximation for this
particular case. In the same figure, a comparison between using the
correction to the panel strength at the end points according to
\eqref{eq:gammacorr}, or by using \eqref{eq:gammaactual} with series
expansions of the correction to make the source distribution really
follow the equation (red line). According to this case, the accuracy
from using the series expansions is actually decreased, showing that
corrections to the end points are a valid approximation, and is hence
the recommended implementation as it is numerically more stable and
faster, as the order of the strength polynomial is kept low and panels
do not have to be split.

\subsection{Computational complexity}

We have seen that under certain conditions, the accuracy can be
increased by using higher order panels. However, the evaluation speed
also needs to be considered before choosing the order of the
panel. Here we will assume that the panels is to be used in a vortex
method, which means that it will be the panel-vortex interaction that
will be the most time-consuming part, since this has to be performed
every time-step, while the BEM matrices often only have to be
calculated once at the start of the simulation.

\subsubsection{Near field evaluation}

\begin{figure}
  \includegraphics[width=1\columnwidth]{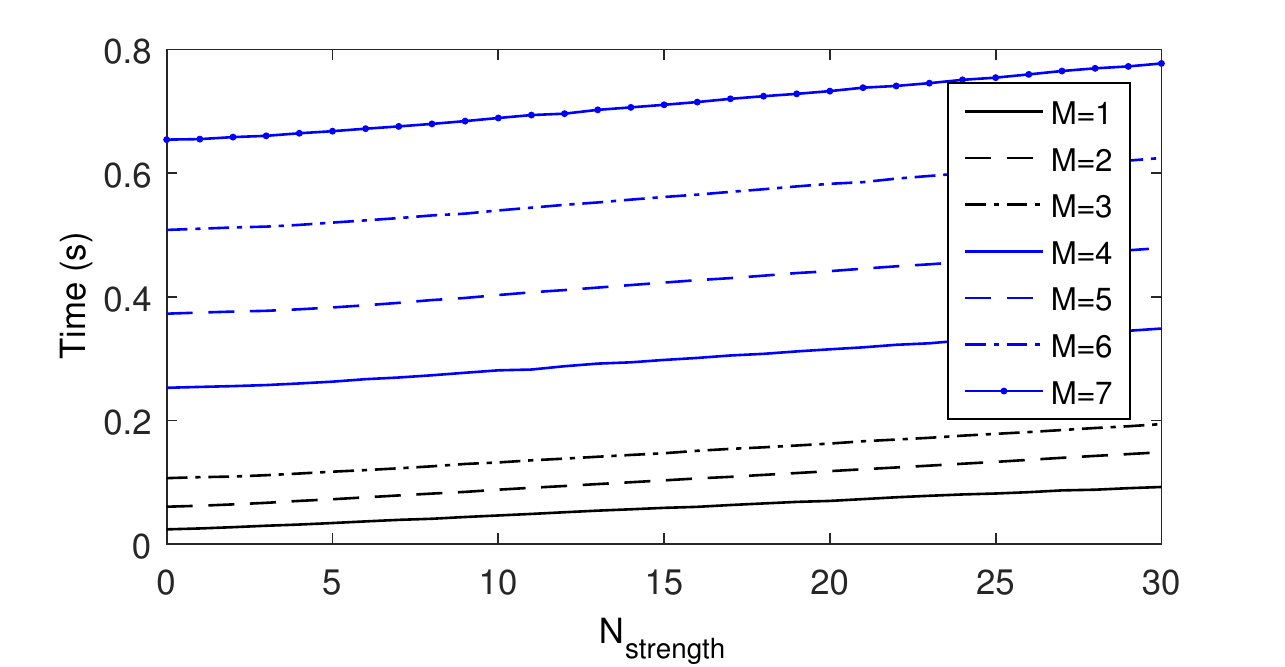}
  \caption{\label{fig:NearFieldTime}Evaluation time for near field
    evaluations. Each line represents a different order for the shape
    of the panel, while different orders of the strength are given on
    the $x$-axis. It can clearly be seen that the order of the shape
    has a much larger impact on the evaluation time than the order of
    the panel strength.}
\end{figure}

In the first case, the evaluation within the near field of the panel
will be considered. Here, 100 panels extending over the diagonal of
the unit square are used, and 1 000 evaluation points are distributed
randomly within the unit square. The time values measured are the wall
clock times for all steps of the calculations; see
Figure~\ref{fig:NearFieldTime} for the results. Clearly, it is the
order of the shape that dominates the time of the evaluation. For
orders above 3, where numerical methods are required to calculate the
roots in \eqref{eq:factorization}, there is a relatively rapid
increase in the evaluation time. From a theoretical point of view, one
can expect the growth both to have a linear component, see
\eqref{eq:generalsolution}, and also a quadratic component from the
root finding algorithm.

\subsubsection{Fast multipole evaluation}

\begin{figure}
  \includegraphics[width=1\columnwidth]{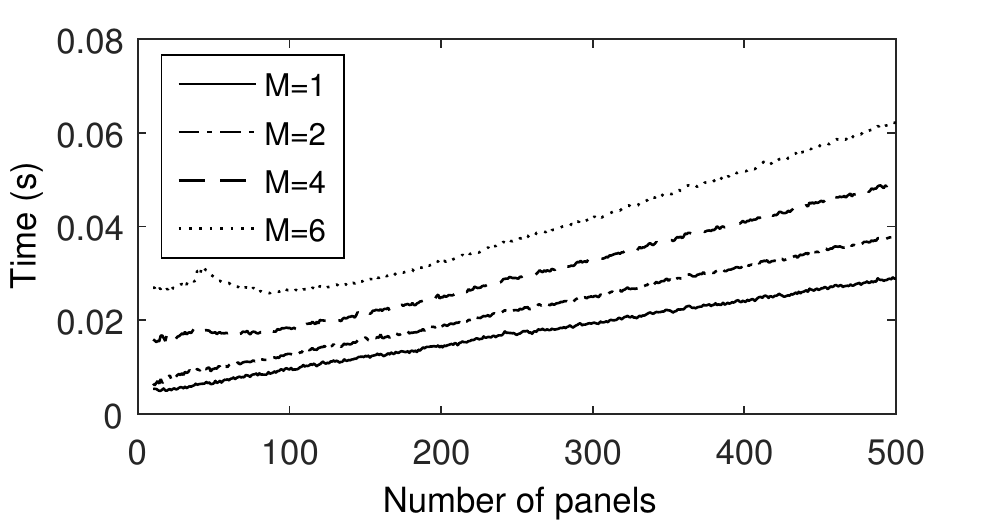}
  \caption{\label{fig:FarFieldTime}Evaluation time with fast multipole
    evaluation of 100 000 particles for the circle in
    \S\ref{sub:unitcircle}. Cubic polynomials for the shape is chosen
    in all cases. The lines represent the added time for including the
    panels into the evaluation. The base time for evaluation without
    the panels is 0.29s. Note that due to symmetry, the highest order
    coefficient for the panel shape will be zero. Hence the values
    $M=2$, 4 and 6, while the corresponding formal accuracy orders are
    3, 5 and 7.}
\end{figure}

In the second case, which is intended to illustrate applications in a
vortex method, the unit circle from \S\ref{sub:unitcircle} will be
used. Here, 100 000 vortices are randomly distributed using a
homogeneous distribution in the polar plane between
$0<\theta\le2\pi,\,1<r\le10$. The fast multipole method is configured
to have about 35 particles in the lowest level boxes, which according
to \cite{Goude13a} is a suitable configuration for the current
implementation for vortex dominated systems. The tolerance of the FMM
is set to $10^{-9}$. In a real vortex method application, the velocity
from the mutual interaction between the vortices also has to be
calculated. This should be carried out in the same FMM evaluation as
the panel interaction to reduce computational overhead. Hence,
evaluation times with the panels included are compared with evaluation
times with only the vortices present. This evaluation method was
chosen to illustrate the added time of using the different panels and
to exclude the time to build the FMM tree of the vortices (which
otherwise will dominate over the vortex evaluation time). As
Figure~\ref{fig:NearFieldTime} showed that it is mainly the panel
shape that determines the evaluation speed, the order of the panel
strength is fixed to 3 in this evaluation. The results are presented
in Figure~\ref{fig:FarFieldTime} and shows that the time indeed
increases with increasing panel order. However, by comparing with the
particle to particle interaction time of 0.29s, the added time of the
panels is still a small part of the total time.

Note that when comparing the evaluation times with
Figure~\ref{fig:NearFieldTime}, 100 times more vortices are included
here, implying about three orders of magnitude faster
evaluation. Without the use of the FMM, the time would be about 0.81s
(100 panels) for all panel orders, as the time for evaluating the
far-field is independent of the panel order.

One trend seen for $M=6$ is that the evaluation time stops to decrease
for low numbers of panels. This is a consequence of both that large
panels give larger regions, where direct evaluation is necessary, and
that large panels will cover several FMM boxes, and hence have to be
split (this is likely the reason for the peak at about 50 panels).


\subsection{Discussion}

One of the main questions is which order of the panels should be used
and as seen from the results this depends on the application. For the
circle case, higher order panels generally give better results, but
for the airfoil, the most suitable panels depend on what the code has
to be used for. If only the circulation is of interest, it can be
suitable to choose a low order panel method that is fast, but higher
order panels have benefits for evaluations close to the airfoil
surface. It is recommended to perform studies of the particular
application before choosing panel method.

One other aspect of the choice between many panels of low order, or
few of high order, is what happens for a real vortex application,
where there are lots of vortices close to the boundary. Here, it can
be beneficial to have many panels, as they more easily can adapt the
strength to fulfill the no-penetration boundary condition. It is seen
in Figure~\ref{fig:FarFieldTime} that for high order panels, using too
few panels will not decrease the evaluation time further. Hence, it
may be a good idea not to use too high panel orders with very few
panels, especially since the high order panels may experience
difficulties when the solution has large derivatives, as is seen at
the trailing edge of the airfoil case.



One additional question, not addressed in this work, is if one should
use source or vortex panels, or a combination of both of them. From an
evaluation point of view, as the whole algorithm is calculated with
complex numbers, using both source and vortex strength should not give
any noticeable increase in computational speed for the velocity
evaluation, however the amount of unknowns will increase. This is an
area where further studies can be performed for best performance.



\section{Conclusions}
\label{sec:conclusions}

We have presented a general high order panel method where, except for
finding the roots of a polynomial, all parts have analytic
solutions. The method hence allows for the evaluation of the panel
contribution with high accuracy for all points in space. We have also
presented how the flow through panels of high order can be calculated,
both from point sources and from other high order panels. The
necessary equations for implementing the panels in a fast multipole
method has been presented as well. Together, these should provide the
necessary basis for including higher order panels in vortex method
applications.

We have demonstrated that for certain applications, a high order
convergence can be obtained with higher order panels. However, in
certain cases, especially close to points with large derivatives of
the solution, numerical difficulties may reduce this high order
convergence, showing that care has to be taken in these cases.

We have also demonstrated that the evaluation time mainly depend on
the order of the shape of the panels, and high order panels are
significantly more expensive to evaluate in the near field. However,
the cost of evaluating the contribution from the higher order panels
quickly decreases for evaluation further away when the fast multipole
method can be applied. If the number of vortices is large, it is
likely that the vortex to vortex interaction will dominate over the
panel to vortex interaction in terms of computational cost.

\subsection{Reproducibility}
\label{subsec:reproducibility}

The full source code of the fast multipole high order panel solver is
available as open source, and comes with a Matlab interface. All test
codes used to generate the results in this work is also included, thus
allowing for easy reproduction of the results. The code is available
at \url{www.stenglib.org}.


\section*{Acknowledgment}

This work was financially supported by the Swedish Energy Agency and
was conducted within StandUp for Wind, a part of the StandUp for
Energy strategic research framework (A.~Goude), and by the Swedish
Research Council within the UPMARC Linnaeus center of Excellence
(S.~Engblom).


\newcommand{\doi}[1]{\href{http://dx.doi.org/#1}{doi:#1}}
\bibliographystyle{unsrt}
\bibliography{refs}

\appendix


\section{Numerical considerations}
\label{app:numerical}

There are many numerical considerations that has to be taken into
account for accurate use of \eqref{eq:generalsolution}. This section
will briefly go through the most important numerical special cases
that have to be handled.

\subsection{Large roots}
\label{subsec:large_roots}

One case that can cause numerical issues is if one or several of the
roots of the polynomial in \eqref{eq:factorization} is large
($\left|x_{k}\right|\gg\lambda$) and this becomes particularly
significant in the case when $N\ge M$. The problem can be
illustrated by looking at the expression for polynomial division
\begin{align}
\frac{\sum_{j=0}^{N}b_{j}\zeta^{j}}{\zeta-R}=\sum_{j=0}^{N-1}d_{j}\zeta^{j}+D,
\end{align}
which has the recursive solution
\begin{eqnarray}
d_{n-1} & = & b_{n}, \\
d_{j} & = & b_{j+1}+R\cdot d_{j+1}.
\end{eqnarray}
Hence the coefficients grow as $R^{j}$ which for large $R$ will cause
large numerical errors (and increasing with the order of the strength
$N$).

To reduce this problem one uses Taylor expansions for the large
roots. Assume that $x_{i}$ is the large root and write
\begin{align}
  \frac{\sum_{j=0}^{N}b_{j}\zeta^{j}}{\prod_{k=1,k\ne i}^{M}\left(\zeta-x_{k}\right)}\frac{1}{\left(\zeta-x_{i}\right)}
 &= -\frac{1}{x_{i}}\frac{\sum_{j=0}^{N}b_{j}\zeta^{j}}{\prod_{k=1,k\ne i}^{M}\left(\zeta-x_{k}\right)}\sum_{m=0}^{\infty}\left(\frac{\zeta}{x_{i}}\right)^{m}.
\end{align}
This expression can be approached in the same way as \eqref{eq:velc},
but now without the issue with the large root.

\subsection{Double roots}

Another case which may produce numerical errors is when two roots are
close to each other, causing division by a small inaccurate number in
\eqref{eq:partialfraccoeff}. This can be handled by studying these two
roots together. Assume that roots $j$ and $k$ are nearly the
same. Then we have
\begin{align}
\nonumber
 A_{k}\log\left(\frac{x_{k}-\lambda}{x_{k}}\right)&+A_{j}\log\left(\frac{x_{j}-\lambda}{x_{j}}\right)
 = \\
 \nonumber
 &\frac{\sum_{m=0}^{M-1}h_{m}x_{k}^{m}}{\left(x_{k}-x_{j}\right)\prod_{n\ne k,j}\left(x_{k}-x_{n}\right)}\log\left(\frac{x_{k}-\lambda}{x_{k}}\right)+ \\
 &\frac{\sum_{m=0}^{M-1}h_{m}x_{j}^{m}}{\left(x_{j}-x_{k}\right)\prod_{n\ne k,j}\left(x_{j}-x_{n}\right)}\log\left(\frac{x_{j}-\lambda}{x_{j}}\right).
\end{align}
By expanding the logarithm of one root around the other we get
\begin{align}
\log\left(\frac{x_{j}-\lambda}{x_{j}}\right)=\log\left(\frac{x_{k}-\lambda}{x_{k}}\right)+\left(x_{j}-x_{k}\right)B_{j},
\intertext{with}
\label{eq:Bj}
B_{j}=\sum_{n=1}^{\infty}\frac{1}{n}\left(x_{j}-x_{k}\right)^{n-1}\left(\left(\frac{1}{x_{k}-\lambda}\right)^{n}-\left(\frac{1}{x_{k}}\right)^{n}\right).
\end{align}
Also, the sum can be expanded as 
\begin{align}
\frac{1}{\prod_{n\ne k,j}\left(x_{j}-x_{n}\right)}&=\frac{1}{\prod_{n\ne k,j}\left(x_{k}-x_{n}\right)}\left(1+\left(x_{k}-x_{j}\right)C_{k,j}\right),
\intertext{with}
C_{k,j}&=\sum_{m=0}^{\infty}\frac{\left(x_{k}-x_{j}\right)^{m-1}}{\left(x_{k}-x_{n}\right)^{m}}.
\end{align}
This allows us to write the whole expansion as
\begin{align}
\nonumber
 & A_{k}\log\left(\frac{x_{k}-\lambda}{x_{k}}\right)+A_{j}\log\left(\frac{x_{j}-\lambda}{x_{j}}\right)\\
 \nonumber
= & \frac{\sum_{m=1}^{M-1}h_{m}\sum_{p=0}^{m-1}x_{k}^{m-1-p}x_{j}^{p}-C_{k,j}\sum_{m=0}^{M-1}h_{m}x_{j}^{m}}{\prod_{n\ne k,j}\left(x_{k}-x_{n}\right)}\log\left(\frac{x_{k}-\lambda}{x_{k}}\right)\\
 & +\frac{\sum_{m=0}^{M-1}h_{m}x_{j}^{m}}{\prod_{n\ne k,j}\left(x_{j}-x_{n}\right)}B_{j}.
\end{align}
For convergence of $B_j$ in \eqref{eq:Bj}, a sufficient condition is
that the distance between the roots is smaller than the distance for
any of the roots to the panel end points.

\subsection{Evaluation close to the panel end}

For connected panels with continuous derivative in the shape and
continuous strength, the velocity solution should be smooth on the
panel surface when shifting between panels. However,
\eqref{eq:generalsolution} contains a weak logarithmic singularity at
this point (when $x_k = 0$ or $\lambda$). Although it is possible to
derive similar techniques as for the double root in this case, as long
as the overall method avoids evaluations in the nearest vicinity of
the end points, this does not appear to be critical within the target
tolerances considered here.


\section{Logarithmic potential}
\label{app:logpot}
The logarithmic potential is generally given by
\begin{equation}
F=\frac{Q+i\Gamma}{2\pi}\log\left(z-z_{v}\right).
\end{equation}
Using the same procedure as for the potential in
\eqref{eq:vortexvelocity} we get
\begin{align}
  \nonumber
  &\frac{1}{2\pi}\int\limits _{0}^{\lambda}\sum_{j=0}^{N}b_{j}\zeta^{j}\log\left(\sum_{k=0}^{M}c_{k}\zeta^{k}\right)d\zeta \\
  \label{eq:logpotfirst}
  &= \frac{1}{2\pi}\left( \int\limits _{0}^{\lambda}\sum_{j=0}^{N}b_{j}\zeta^{j}\sum_{k=1}^{M}\log\left(\zeta-x_{k}\right)d\zeta-\sum_{j=0}^{N}\frac{b_{j}}{j+1}\lambda^{j+1}\log\left(c_{M}\right)\right) +N_{2}i\int\limits _{0}^{\lambda}q
\end{align}
In terms of
\begin{align}
\int\limits _{0}^{\lambda}\sum_{j=0}^{N}b_{j}\zeta^{j}\log\left(\zeta-x_{k}\right)d\zeta &= \left[\sigma_{0}\left(\zeta;x_{k}\right)\log\left(\zeta-x_{k}\right)\right]_{0}^{\lambda}-\sum_{m=1}^{N+1}\sigma_{m}\left(x_{k}\right)\lambda^{m},
\label{eq:logpotmiddle}
\end{align}
where
\begin{align*}
\sigma_{0}\left(\zeta;x_{k}\right) &:= \sum_{j=0}^{N}\frac{b_{j}}{j+1}\left(\zeta^{j+1}-x_{k}^{j+1}\right),\\
\sigma_{m}\left(x_{k}\right) &:= \frac{1}{m}\sum_{j=m-1}^{N+1}b_{j}\frac{1}{\left(j+1\right)}x_{k}^{j+1-m}.
\end{align*}
Combining \eqref{eq:logpotfirst} with \eqref{eq:logpotmiddle} we get
\begin{align}
 \nonumber
 \int\limits _{0}^{\lambda}\sum_{j=0}^{N}b_{j}\zeta^{j}&\log\left(\sum_{k=0}^{M}c_{k}\zeta^{k}\right)d\zeta = \sum_{j=0}^{N}\frac{b_{j}}{j+1}\left(\lambda^{j+1}\right)\log\left(\sum_{k=0}^{M}c_{k}\lambda^{k}\right) \\
 \label{eq:logpot}
 &-\sum_{k=1}^{M}\left(\sum_{j=0}^{N}\frac{b_{j}}{j+1}\left(x_{k}^{j+1}\right)\log\left(\frac{x_{k}-\lambda}{x_{k}}\right)-\sum_{m=1}^{N}\sigma_{m}\left(x_{k}\right)\lambda^{m}\right),
\end{align}
which is the final expression for the logarithmic potential. As in
\S\ref{subsec:large_roots}, in case of large roots, a series expansion
might be employed for numerical reasons.


\section{Convergence with respect to number of panels}
\label{app:convergence}

This section will cover the convergence estimate for the far-field
region. The flow velocity of the panel is an expansion on the form
\begin{equation}
  \label{eq:VC}
\overline{V}=\int\limits _{0}^{\lambda}\frac{\sum_{k=0}^{N}b_{k}\zeta^{k}+O\left(\zeta^{N+1}\right)}{z'-\zeta-i\sum_{k=0}^{M}a_{k}\zeta^{k}+O\left(\zeta^{M+1}\right)}d\zeta.
\end{equation}
In the far-field, where $z'$ is large enough, one can readily expand
the denominator in a geometric series containing negative powers of
$z'$. Hence by the truncation of polynomials in both the nominator and
the denominator, one finds that the far-field convergence rate is
$\min\left(M+1,N+1\right)$.

However, in terms of the strength, a higher convergence rate can be
obtained than the suggested value $N+1$. In contrast to the shape, the
source strength is not obtained from an explicit value. Instead, it is
calculated implicitly by inverting the BEM matrix. For the far-field,
we can look at this problem in terms of series expansions. Write
\eqref{eq:VC} as an expansion around some point $z_0$,
\begin{equation}
  \overline{V\left(z\right)}=
  \sum_{n=1}^{\infty}\frac{f_{n}}{\left(z-z_{0}\right)^{n}}.
\end{equation}
Similarly, write the expansion of the true source as
\begin{equation}
  \label{eq:convergenceG}
  G\left(\zeta\right)=\sum_{k=0}^{\infty}g_{k}\zeta^{k}.
\end{equation}
It follows that the best possible coefficients, in the sense of the
highest possible truncation order, are given by
(cf.~\eqref{eq:fnseries})
\begin{equation}
  \label{eq:fnC}
  f_{n}=\int\limits _{0}^{\lambda}G\left(\zeta\right)\left(\left(\zeta+i\sum_{k=0}^{M}a_{k}\zeta^{k}\right)e^{i\theta}-z_{0}+z_{1}\right)^{n-1}d\zeta.
\end{equation}
Note that the constant terms $z_{0}$, $z_{1}$ and $a_0$ in
\eqref{eq:fnC} represent the distance between the starting point and
the center of the panel, and will be of the order $\lambda/2$.

Combining \eqref{eq:convergenceG} and \eqref{eq:fnC}, we can write the
values for the coefficients as
\begin{equation}
  \label{eq:fnnC}
  f_{n}=\sum_{k=0}^{\infty}g_{k}\lambda^{k+n}h_{n,k}(\lambda),
\end{equation}
where $h_{n,k}(\lambda) = O(1)+O(\lambda)$. Note that for panels with
linear shape, $h_{n,k}$ would be independent of $\lambda$. Now, assume
that we want to approximate the source with a polynomial
\begin{equation}
  G(\zeta) \approx B\left(\zeta\right)=\sum_{k=0}^{N}b_{k}\zeta^{k}.
\end{equation}
Instead of truncating $G$ in \eqref{eq:convergenceG} to $N$
coefficients, we will instead choose the coefficients of $B$ to
exactly match the $N+1$ lowest terms in the series expansion
\eqref{eq:fnnC}. This gives the system of equations {\scriptsize
\begin{equation}
\label{eq:exact_terms}
\left(\begin{array}{cccc}
\lambda h_{1,0} & \lambda^{2}h_{1,1} & \dots & \lambda^{N+1}h_{1,N}\\
\lambda^{2}h_{2,0} & \lambda^{3}h_{2,1} & \dots & \lambda^{N+2}h_{2,N}\\
\vdots & \vdots & \ddots & \vdots\\
\lambda^{N+1}h_{N,0} & \lambda^{N+2}h_{N,1} & \dots & \lambda^{2N+1}h_{N,N}
\end{array}\right)\left(\begin{array}{c}
b_{0}\\
b_{1}\\
\vdots\\
b_{N}
\end{array}\right)=\left(\begin{array}{c}
\sum_{k=0}^{\infty}g_{k}\lambda^{k+1}h_{1,k}\\
\sum_{k=0}^{\infty}g_{k}\lambda^{k+2}h_{2,k}\\
\vdots\\
\sum_{k=0}^{\infty}g_{k}\lambda^{k+N+1}h_{N,k}
\end{array}\right),
\end{equation}}
which we can write in matrix notation as
\begin{equation}
  Ab=c \Longrightarrow b=A^{-1}c.
\end{equation}
We can now split the vector $c$ into two parts
\begin{equation}
  c_{m}=c_{1,m}+c_{2,m}=\sum_{k=0}^{N}g_{k}\lambda^{k+m+1}h_{m+1,k}+\sum_{k=N+1}^{\infty}g_{k}\lambda^{k+m+1}h_{m+1,k},
\end{equation}
where $m$ is the element index, starting from 0. Noting that $c_1 =
Ag$ we get
\begin{equation}
\label{eq:Csplit}
b=A^{-1}\left(c_{1}+c_{2}\right)=g+A^{-1}c_{2} =: g+b_2,
\end{equation}
where $b_2$ can be thought of as the error coefficients. As the lower
order terms up to $\lambda^{N+1}$ belongs to $c_1$, while $b_2$
contains the remaining terms from $c_2$, the lowest order term in
$b_2$ will be of order $b_{2,0} = O(\lambda^{N+1})$. Now, for any
coefficient $f_{n}$, we have the approximation
\begin{equation}
  f_n \approx 
  \hat{f}_{n}=\sum_{k=0}^{N}b_{k}\lambda^{k+n}h_{n,k}=\sum_{k=0}^{N}g_{k}\lambda^{k+n}h_{n,k}+\sum_{k=0}^{N}b_{2,k}\lambda^{k+n}h_{n,k}.
\end{equation}
The error is just the difference between the approximation and the
correct solution, thus,
\begin{align}
  \label{eq:convergence_error}
  \hat{f}_{n}-f_{n} &= 
  \sum_{k=0}^{N}b_{2,k}\lambda^{k+n}h_{n,k}-\sum_{k=N+1}^{\infty}g_{k}\lambda^{k+n}h_{n,k}.
\end{align}
The lower order terms up to $n=N+1$ in \eqref{eq:convergence_error}
are 0 thanks to the matching in \eqref{eq:exact_terms}. Also, since
the lowest order term in $b_{2}$ is $O(\lambda^{N+1})$, the lowest
order term in \eqref{eq:convergence_error} is $O(\lambda^{n+N+1})$. It
follows that the lowest order term we can find in any of the
coefficients is of order $\lambda^{2N+3}$. Finally, one notes that the
number of panels increases with $1/\lambda$, and hence the total
convergence rate in the far-field for the source strength is bounded
by $\lambda^{2N+2}$.

\end{document}